\documentclass[12pt]{article}
\usepackage{amssymb}

\newcommand{\f}{\varphi}

\newcommand{\pp}{\mathbb{P}}%{{\bf P}}
\newcommand{\zz}{\mathbb{Z}}%{{\bf Z}}

\newcommand{\aA}{{\cal A}}
\newcommand{\bB}{{\cal B}}
\newcommand{\cC}{{\cal C}}
\newcommand{\dD}{{\cal D}}
\newcommand{\oO}{{\cal O}}
\newcommand{\tT}{{\cal T}}

\newcommand{\lle}{\mbox{\raisebox{0.25ex}{${\scriptscriptstyle\le}$}}}
\newcommand{\gge}{\mbox{\raisebox{0.25ex}{${\scriptscriptstyle\ge}$}}}
\newcommand{\Dlz}{{\cal D}^{\lle0}}
\newcommand{\Dgo}{{\cal D}^{\gge1}}
\newcommand{\tr}{{\mbox{$t$-struc}\-tur}}
\newcommand{\Tn}{{\mbox{$T$-in}\-va\-ri\-an}}
\newcommand{\bcdot}{{\mbox{\boldmath{$\cdot$}}}}

\newcommand{\Hom}{\mathop{\rm Hom}\nolimits}
\newcommand{\Ext}{\mathop{\rm Ext}\nolimits}

\newcommand{\rk}{\mathop{\rm rk}\nolimits}

\newcommand{\naive}{naive}

\newcommand{\simto}{\mathbin{\lefteqn{%
\mbox{\raisebox{0.6ex}{$\hspace{1pt}\sim$}}}%
\mbox{\raisebox{-0.4ex}{$\to$}}}}

\newtheorem{LEM}{Lemma}
\newtheorem{REM}[LEM]{Remark}
\newtheorem{THM}[LEM]{Theorem}
\newtheorem{PROP}[LEM]{Proposition}

\newenvironment{ITEM}[1]{\begin{trivlist}\item[\hspace{\labelsep}%
{\bf#1}]}{\end{trivlist}}

\newenvironment{PROOF}[1]{\par\noindent{\bf Proof#1}}{\par\medskip}

\begin{document}

%\tableofcontents  % <--- {\aa}{\'O}{\`I}{\'E} {\^I}{\~O}{\"O}{\^I}{\"I} {\'O}{\"I}{\"A}{\AA}{\`O}{\"O}{\'A}{\^I}{\'E}{\AA}, {\'O}{\"I}{\^O}{\`O}{\'E} "%" + {\^I}{\'A}Ё{\'A}{\`I}{\AA} {\'O}{\^O}{\`O}{\"I}{\"E}{\'E}.

\title{Operations on \tr es and perverse coherent sheaves}
\author{Alexey Bondal\\
%{\small Steklov Mathematical Institute}\\[-1ex]
%{\small Gubkina 8, GSP1}\\[-1ex]
%{\small 117966 Moscow}\\[-1ex]
%{\small Russia}\\[-1ex]
%{\small {\tt e-mail: bondal@mi.ras.ru}}
}
\date{}

\maketitle

{\bf {\Large Introduction.}}

%+++++++++++++++++++++++++++++++++++++++++++++++++++++++++++++++++++++++++++

Suppose we have an equivalence $\dD ({\cal A})\simeq \dD ({\cal
B})$ between the derived categories of two abelian categories
$\aA$ and $\bB$. How this helps to study $\dD ({\cal A})$?
The virtue is that we can transport
the standard \tr e from $\dD ({\cal B})$ into $\dD ({\cal A})$
via the equivalence and look how it interplays with
the standard \tr e in $\dD ({\cal A})$. In particular, we can try to produce
new \tr es from these two.

In this paper we consider two binary operations on \tr es in a
triangulated category, $ \dD$. More precisely, there is a natural
partial order on the set of \tr es and the operations are just
that of intersection and union as they are defined in the lattice
theory. The basic problem is to prove the existence of
intersection and union for two given \tr es.

We define {\em lower and upper consistent pairs} of \tr es. These
are some elementary conditions which guarantee the existence
of intersection or union for the pairs. Then we take a triple of \tr es
and examine various conditions of consistency in pairs of \tr es
in this triple which imply consistency in the pair that comprises
one of them and the intersection (or union) \tr e of the other
two. This gives some rules for iterating the operations of union and intersection.

The lattice of vector subspaces in a vector space does not satisfy distributivity law, but rather
{\em modularity law}. The partially ordered set of \tr es is somewhat similar. We prove, in particular, some versions of modularity law under suitable
conditions on consistency in pairs. We consider abstract partially ordered sets with two binary relations subject to  axioms, which upper and low consistency obey. We call them {\em sets with consistencies}. We prove that such a set $P$ allows a map into a universal {\em lattice} with consistencies, $U(P)$.

A {\em chain} in a poset is a fully ordered set of elements.
A theorem due to Birkhoff \cite{Bir} claims that two chains in a modular lattice generate a distributive lattice.
A pair of chains $\{a_i\}$ and $\{b_j\}$ is said to be consistent if all pairs $(a_i, b_j)$ are both upper and low consistent. We prove a theorem that a consistent pair of chains in a set with consistencies generates a distributive lattice $L$.

Our application of the developed techniques is to perverse coherent sheaves on schemes of finite type over a field.
We consider the derived category $\dD (X)=\dD^b_{\rm coh}(X)$ of coherent sheaves on such a scheme $X$. It has the standard t-structure. The category possesses Grothendieck-Serre duality functor $D: \dD (X)\to \dD (X)^{opp}$, defined for ${\cal F}\in \dD (X)$ by
$$
D{\cal F}={\cal H}om_X({\cal F}, \omega^{\bcdot}_X),
$$
where $\omega^{\bcdot}_X$ stands for the dualizing complex on $X$.
Since $D$ is an anti-equivalence, it allows to define another, `dual', t-structure on $\dD (X)$ by transporting the standard t-structure via $D$. We consider two chains of t-structures: one is obtained by shifts of the standard t-structure and the other one by shifts of the dual t-structure. We prove that this pair of chains is consistent. Therefore, they generate a distributive lattice of t-structures by the above version of Birkhoff theorem. It is a set of perverse t-structures on coherent sheaves (those which are independent of stratification of $X$) considered independently by Deligne, Bezrukavnikov \cite{Bez}, Kashiwara \cite{Kash}, Gabber \cite{Gab}.

The hard part in constructing a \tr e with required properties is to prove existence of the
adjoint to the embedding of what-is-to-be $\Dlz$ into $\dD$. The standard way to
produce such a functor is by using limits in the category.
%Here we should outline what type of categories we are interested in. Let us mention only one motivating problem.

According to the proposal in paper  \cite{BO}, derived categories of coherent
sheaves can play a crucial role in the Minimal Model Program of
birational algebraic geometry.
%(cf. \cite{BO} and reference therein).
We expect that important ingredients of the Minimal Model Program, such
as flips and flops, could be constructed and understood  by means of transformation
of \tr es (see \cite{Br1} for 3-dimensional flops).
%The landscape
%should be enriched by involving noncommutative varieties (cf.
%\cite{BO},\cite{VdB1}, \cite{VdB2}).

This problem requires considering triangulated categories of `small' size,
like those arising in `realistic' algebraic geometry, such as the
bounded derived categories of complexes of coherent sheaves on
algebraic varieties as opposed to the unbounded derived categories of
quasi-coherent sheaves. There is no suitable general properties on
existence of limits in such categories. For that reason, standard
constructions often lead us {\em a priori} beyond the `small'
category. For example, the derived functor of cohomology with
support takes, in general, coherent complexes to quasi-coherent
ones.

One of the main technical problems in studying algebraic geometry by
means of triangulated categories can be formulated as to find appropriate finiteness properties on
suitable subcategories in the derived categories of coherent sheaves and keep track of these properties
under relevant transformations of these subcategories.

%Examples of such
%properties for abelian categories are noetherianess and
%artinianess of cat. To mention a one for triangulated categories, there
%is a notion of strongly finitely generated triangulated categories
%\cite{BVdB}, which might be considered as a basic property of the
%derived category of perfect complexes on a smooth (noncommutative)
%algebraic variety.
We hope that this paper shed light on the interplay of \tr es, thus giving approach to understand the finiteness conditions relevant to \tr es.

It would be interesting to compare the approach of this paper with Tom Bridgeland's work on producing stability conditions (hence a plenty of \tr es!) via deformation argument \cite{Br2}.

I am indebt to Michel~Van-den-Bergh for useful discussions.

This work was done more than 15 years ago and was, for the first time, reported on Shafarevich's seminar at Steklov Institute in late 90's. It is my pleasure to devote this paper to I.R. Shafarevich on the occasion of his 90-th anniversary. The influence of Shafarevich on Russian algebraic geometry can hardly be overestimated.

%---------------------------------------------------------------------------
\section{Consistent systems of \tr es}
\subsection{The partially ordered set of \tr es}
\label{poset}
Let $k$ be any field, $\dD$ a $k$-linear triangulated category with an
auto-equivalence $T$, called translation (or shift) functor (see~\cite{V}).

For a subcategory $\cC$ in $\dD$ its
{\em right} (resp. {\em left}) {\em orthogonal} $\cC^\bot$ (${}^\bot\cC$) is
the strictly full subcategory in $\dD$ of objects $X$ such
that $\Hom(\cC,X)=0$ (resp. $\Hom(X,\cC)=0$).

We use the standard notation $X[k]:=T^kX$, ${k\in\zz}$.

Recall (see~\cite{BBD}) that a {\it\tr e} in a triangulated category~$\dD$ is a pair of strictly
full subcategories $(\Dlz,\Dgo)$, satisfying the following conditions:
\begin{itemize}
\item[{\it i})] $T\Dlz\subset\Dlz$, $T^{-1}\Dgo\subset\Dgo$;
\item[{\it ii})] $\Hom(\Dlz,\Dgo)=0$;
\item[{\it iii})] for any ${X\in\dD}$ there exist $X_{-}\in\Dlz$, $X_{+}\in\Dgo$
and an exact triangle
\begin{equation}
X_{-}\to X\to X_{+}\to X_{-}[1]
\label{f1}
\end{equation}

\end{itemize}

The basic example of a \tr e comes up when $\dD$ is equivalent to the bounded derived category $\dD^b(\aA)$ of
some abelian category~$\aA$. For this case $\Dgo$ (resp. $\Dlz$) comprises the complexes
with cohomology in positive (resp. nonpositive) degree.

Categorically, the main virtue of the notion of \tr e is that in
view of conditions (i) and~(ii) the triangle~(\ref{f1}) for a
given~$X$ is unique up to  unique isomorphism~\cite{BBD}. In other
words, formulas
$$
\begin{array}{l}
\tau_{\lle0}(X)=X_{-},\ \
\tau_{\gge1}(X)=X_{+}
\end{array}
$$
define the truncation functors $\tau_{\lle0}\colon{\dD\to\Dlz}$,
$\tau_{\gge1}\colon{\dD\to\Dgo}$, which are easily seen to be
respectively the right and the left adjoint to the inclusion
functors $i_{\lle0}\colon$ ${\Dlz\to\dD}$, $i_{\gge1}\colon$
${\Dgo\to\dD}$. This suggest the following reformulation of the
data defining a \tr e.

\begin{LEM}\label{l1}
To define a \tr e in a triangulated category~$\dD$ is equivalent
to one of the following data:

\begin{description}
\item[Data 1.] A strictly full subcategory $\Dlz$ in $\dD$,
closed under the action of~$T$ and such that the inclusion functor
$i_{\lle0}\colon{\Dlz\to\dD}$ has a right adjoint~$\tau_{\lle0}$,
\item[Data 2.] A strictly full subcategory $\Dgo$ in $\dD$,
closed under the action of~$T^{-1}$ and such that the inclusion
functor $i_{\gge1}\colon{\Dgo\to\dD}$ has a left
adjoint~$\tau_{\gge1}$.
\end{description}
\end{LEM}

\begin{PROOF}{.} As we have already constructed adjoint functors for
 a given \tr e, it remains to explain how to construct a \tr e from say data~1.

Put $\Dgo =(\Dlz )^{\bot}$. As ${T\Dlz\subset\Dlz}$,
then obviously ${T^{-1}\Dgo\subset\Dgo}$.

For any
${X\in\dD}$, the adjunction morphism $i_{\lle0}\tau_{\lle0}{X\to X}$ can be inserted in a triangle:
$$
i_{\lle0}\tau_{\lle0}X\to X\to X_{+}\to i_{\lle0}\tau_{\lle0}X[1].
$$
$X_{+}$ is easily seen to be in~$\Dgo$.
Hence the triangle is of the form~(\ref{f1}).
\end{PROOF}

This simple observation is useful in constructing  \tr es.

\begin{REM}\label{1rem}\rm
Note that in general neither inclusion nor truncation functors are
exact. But one can easily see that if $A\to{B\to C}$ is an exact
triangle with $A$, ${C\in\Dlz}$ (resp.~$\Dgo$), then ${B\in\Dlz}$
(resp.~$\Dgo$). It follows  by `rotation of triangle' that if $A$
and $B$ are in $\Dlz$ (resp. $B$, ${C\in\Dgo}$), then also
${C\in\Dlz}$ (resp. ${A\in\Dgo}$).
\end{REM}

Any \tr e generates two sequences of subcategories $\dD^{\lle
n}:=\Dlz[-n]$, $\dD^{\gge n}:=\Dgo[1-n]$. The subcategory
$\cC={\Dlz\cap\dD^{\gge0}}$ has a natural structure of  abelian
category (see~\cite{BBD}), it is called the {\it heart of \tr
e.} The objects which lie in the heart (and also their shifts)
are called {\em pure}.

When ${\cal D}={\cal D}^b({\cal A})$, the heart of the standard
\tr e comprises the complexes over ${\cal A}$ with cohomology in degree zero.
Thus, the heart is
equivalent to the primary abelian category ${\cal A}$. We identify
${\cal A}$ with this full subcategory in ${\cal D}^b({\cal A})$.

\vspace{0.5cm}

If $\dD$ is  triangulated, then the opposite category $\dD ^{op}$
is also naturally triangulated. Namely, $T_{\dD
^{op}}=T^{-1}_{\dD}$  and triangles in $\dD ^{op}$ are defined by
inverting arrows in triangles in $\dD$. A $t$-structure $\tT$ in
$\dD$ defines a t-structure $\tT ^{op}$ in $\dD ^{op}$ by
exchanging the roles of $\Dlz$ and $\Dgo$. The heart of $\tT^{op}$
is the opposite abelian category to the heart of $\tT$ .

\vspace{0.5cm}

An important special case is a \tr e satisfying $T^{-1}\Dlz \in
\Dlz$. It follows from the remark~\ref{1rem} that for such a \tr e
both $\Dlz$ and $\Dgo$ are triangulated subcategories. We call
such a \tr e {\it translation invariant\/} or, simply, {\it \Tn
t\/}.

A strictly full triangulated subcategory $\aA$ in~$\dD$ which has
the right (resp. left) adjoint to the inclusion
$i\colon{\aA\to\dD}$ is called  {\it right\/} (resp. {\it left})
{\it admissible\/} ~\cite{B}. For a \Tn t \tr e, $\Dlz$
(resp.~$\Dgo$) is right (resp. left) admissible. This gives a 1-1
correspondence between \Tn t \tr es and right (or left) admissible
subcategories.

For a \Tn t \tr e all the subcategories $\dD^{\lle k}$,
${k\in\zz}$, are equal (the same for $\dD^{\gge k}$, ${k\in\zz}$).
Hence the heart of a \Tn t \tr e is zero.

\vspace{0.5cm}

Given a triangulated category~$\dD$,
consider the following partial order on the set of its \tr es
(we ignore all possible set theoretical questions which arize from the fact that
the objects of a category comprise a class not a set,
as being irrelevant to our discussion).

For a pair $(\tT_1,\tT_2)$ of \tr es, we say
${\tT_1\subseteq\tT_2}$ if ${\Dlz_1\subseteq\Dlz_2}$. Then also
 ${\Dgo_2\subseteq\Dgo_1}$.

The minimal element~$\bf 0$ in this ordered set is the \tr e with ${\Dlz=0}$;
the maximal one, $\bf 1$, has ${\Dlz=\dD}$.

Following the standard definition of union and intersection for partially ordered sets we introduce:
\begin{ITEM}{Definition}
The {\it (abstract)
intersection\/}
$\bigcap\limits_{i\in I}\tT_i$ of a set $\{\tT_i\}_{i\in I}$ of \tr es
in~$\dD$ is the supremum of \tr es $\tT$ in~$\dD$ such that
${\tT\subseteq\tT_i}$, for all ${i\in I}$. If it exists, it is unique.
The {\it (abstract) union\/} $\bigcup\limits_{i\in I}\tT_i$ is defined dually.
\end{ITEM}

We also denote   the intersection and union operations by $\cdot$
and~$+$ and call them  {\it product\/} and {\it sum}. Indeed,
these are just product and coproduct in the category related to
the partially ordered set.

The intersection and the union of \Tn t \tr es, when exists, is
\Tn t.

There is a natural `\naive' candidate for $\Dlz_\cap$ (resp.
$\Dgo_\cup$) for the intersection (resp. union) \tr e of a set~$\{\tT_i\}$:
$\Dlz_\cap=\bigcap\limits_{i\in I}\Dlz_i$, (resp.
$\Dgo_\cup=\bigcap\limits_{i\in I}\Dgo_i$).

There is a duality functor on the category of partially ordered sets that inverts
the partial order. It takes the union operation into the intersection and vice versa.
We refer often to this duality in the present paper.

\subsection{An example}
\label{example}

The following example of geometric nature provides with a pair
$(\tT_1,\tT_2)$ of \tr es such that ${\Dlz_1\cap\Dlz_2}$ does not
have right adjoint. This shows that the `\naive' intersection does
not necessarily coincides with the abstract one.

%\begin{ITEM}{Example}
Here we assume $k$ to be algebraically closed and of characteristic zero.
Let $\dD=\dD^b_{\rm coh}(\pp^2)$ be the bounded derived category of coherent
sheaves on the projective plane $\pp^2=\pp(V)$, $\dim{V=3}$. We consider two
translation invariant \tr es in~$\dD$.
Define corresponding right admissible triangulated
subcategories $\Dlz_1$ and~$\Dlz_2$ by

$$
\Dlz_1=\langle\oO(-2),\oO(-1)\rangle,\qquad
\Dlz_2=\langle\oO(1),\oO(2)\rangle,
$$
i.\,e. the minimal strictly full triangulated subcategories generated by
the indicated  sheaves.
% $(\oO(-2),\oO(-1))$ and $(\oO(1),\oO(2))$ respectively.

%%%%%%%%%%%%%%%%%%%%%%%%%%%%%%%%%%%%%%%%%%%%%%%%%

\begin{LEM}\label{lemex}

The category $\dD_\cap=\Dlz_1\cap\Dlz_2$ is not right admissible.

\end{LEM}

\begin{PROOF}{.}
%\end{PROOF}
The category $\Dlz_2$ is equivalent as a triangulated category to
$\dD^b({\bmod}{-}A)$, the bounded derived category of representations of the
path algebra of the quiver (see~\cite{B}):
$$
\bcdot\stackrel{\displaystyle V}{\longrightarrow}\bcdot
$$

The equivalence can be given by the functor that assigns to a complex of
representations
$\stackrel{\lefteqn{\scriptstyle U^{\bcdot}}
 \hphantom{\scriptstyle U}}{\bcdot}
 \longrightarrow
 \stackrel{\lefteqn{\scriptstyle W^{\bcdot}}
 \hphantom{\scriptstyle W}}{\bcdot}$
with the structure map $\f\colon U^{\bcdot}\otimes{V\to
W^{\bcdot}}$, the complex of sheaves
$U^{\bcdot}\otimes\oO(1)\to{W^{\bcdot}\otimes\oO(2)}$ in~$\Dlz_2$
with the differential induced by~$\f$.

Since $A$ is hereditary, any indecomposable object in
$\dD^b({\bmod}{-}A)$ is pure, i.\,e. it is defined by a  map
$\f\colon U\otimes{V\to W}$, with $U$, $W$  vector spaces (up to a
common shift). The corresponding complex of sheaves in~$\Dlz_2$
has the form
\begin{equation}
U\otimes\oO(1)\to W\otimes\oO(2)
\label{fd2}
\end{equation}
with the components in degrees~0 and~1.

Let us describe indecomposables which are in~$\dD_\cap$. Note that
an object ${X\in\dD}$ is in~$\Dlz_1$ iff ${\Ext^\bcdot(\oO,X)=0}$.
Applying this to~(\ref{fd2}) we get that the map
$$
\widehat{\f}\colon U\otimes V^{*}\to W\otimes S^2V^{*},
$$
obtained from $\f$ by partial dualization $\widetilde{\f}\colon
U\to W\otimes V^{*}$, tensoring with the identity map in $V^{*}$
and then by partial symmetrization, has to be an isomorphism. As
${\dim V=3}$, we see that for $X\in\dD_\cap=\Dlz_1\cap\Dlz_2$ we
have $\dim U=2\dim W$.

If
we take $W$ to be 1-di\-men\-sion\-al, then ${\dim U=2}$, hence the image
of~$\widetilde{\f}$ annihilates a line in~$V$. It follows that the image
of~$\widehat{\f}$ consists of quadrics that vanish on this line, i.\,e.
$\widehat{\f}$ is not an isomorphism. Therefore, there is no
indecomposable  in~$\dD_\cap$ with ${\dim W=1}$.

Yet, $\dD_\cap$ is not empty. Indeed, take any exceptional vector
bundle on~$\pp^2$, i.\,e. such that
\begin{equation}
\Hom(E,E)=k,\qquad \Ext^i(E,E)=0,\quad\mbox{for }\ i\ne0.
\label{exep}
\end{equation}

Consider the sheaf $F={\cal E}nd_0E$ of traceless local endomorphisms
of~$E$.
 It
is non-zero, if ${\rk E>1}$. Then (\ref{exep}) implies
$\Ext^\bcdot(\oO,F)={{\rm H}^\bcdot(F)=0}$. Since $\Dlz_1=\oO^{\bot}$
in~$\dD$, then ${F\in\Dlz_1}$. As ${F^*=F}$, then
$\Ext^\bcdot(F,\oO)={{\rm H}^\bcdot(F^*)=0}$. Since
$\Dlz_2={}^{\bot}\oO$, then ${F\in\Dlz_2}$, i.\,e.
${F\in\dD_\cap}$.

%An alternative example of non-zero objects in the category $\dD_\cap$ was proposed by the referee.
%One can take $F$ to be any stable rank 2 vector bundle on $\pp^2$ with Chern classes $c_1(F)=0$ and $c_2(F)=2$.
%These bundles can be constructed by means of {\em generic} extensions of the form:
%$$
%0\to \oO (-1)\to F \to I_Z(1)\to 0
%$$
%where $I_Z$ is the ideal sheaf of the subscheme of a triple of non collinear points in $\pp^2$
%(this is basically Serre construction).
%One can directly see from this short exact sequence that ${\rm H}^\bcdot(F)=0$ is trivial.
%Since $F$ is a rank 2 vector bundle with $c_1(F)=0$, it is self-dual, which implies that ${\rm H}^\bcdot(F^*)=0$,
%hence ${F\in\dD_\cap}$.

Now take $G=\oO(-1)[1]$. If ${F\in\dD_\cap}$ is of the form~(\ref{fd2}),
then
%the standard formulas for the cohomology of line bundles on~$\pp^2$
%shows that:
\begin{equation}
\Hom^1(F,G)=W^*,\qquad \Hom^i(F,G)=0,\quad\mbox{for }\ i\ne1.
\label{homfg}
\end{equation}
Suppose that $\dD_\cap$ is right admissible and $\tau$ is the right adjoint
to the embedding functor ${\dD_\cap\hookrightarrow\dD}$. Put $G'=\tau G$.
Then $G'\in\dD_\cap$. If $G'$ is of the form
$U'\otimes\oO(1)\to{W'\otimes\oO(2)}$, while $F$ is of the form~(\ref{fd2}),
then $\Hom^\bcdot(F,G')$ is the homology of the complex
\begin{equation}
U^*\otimes U'\oplus W^*\otimes W'\to U^*\otimes W'\otimes V^*,
\label{homfp}
\end{equation}
with two non-trivial components in degrees~0 and~1. Since $\dim U=2\dim W$
and $\dim U'=2\dim W'$, the Euler characteristic of this complex is
$-\dim W\cdot\dim W'$.

From the adjunction property and from~(\ref{homfg}) we conclude:
$$
-\dim W=\chi(F,G)=\chi(F,G')=-\dim W\cdot\dim W'.
$$
Therefore, $\dim W'$ must be equal to~1, which is, as we have seen,
impossible.

Generally, $G'$, being an object in~$\dD_\cap$, has a decomposition into
a sum of
objects of the form~(\ref{fd2}) and their translations. Any non-tri\-vial
translation would give a non-tri\-vial contribution in $\Hom^i(F,G')$, for
${i\ne1}$, because these $\Hom$-groups are calculated by translations of the
complex~(\ref{homfp}). But this contradicts the adjunction property for $\tau
$. Therefore, there is no chance for $G'$ to exist. The lemma is proved.
\end{PROOF}

One shows in the same way that the abstract intersection
${\tT_1\cap\tT_2}$ is~$\bf 0$.
%\end{ITEM}

%---------------------------------------------------------------------------
\subsection{Consistent pairs of \tr es and operations}
\label{pairs}

Here we define lower and upper consistent pairs of \tr es. There are various possible ways for generalizing these notion, which we hope to address in an appropriate place.

%For a sequence of \tr es $\tT_i=(\Dlz_i,\Dgo_i)$, ${i\in I}$, numerated by a
%set~$I$, in order to distinguish different
We mark truncation and inclusion functors of \tr es with relevant superscripts.
% (as $\tau^i_{\lle0}$, $i^j_{\gge0}$).

\begin{ITEM}{Definition}
A {\it pair} of \tr es $(\tT_1,\tT_2)$ is called {\it lower
consistent} if
$$
\tau^1_{\lle0}\Dlz_2\subset\Dlz_2,
$$
and {\it upper consistent} if
$$
\tau^2_{\gge1}\Dgo_1\subset\Dgo_1.
$$
It is called {\it consistent} if it is simultaneously lower and upper
consistent.
\end{ITEM}

For consistent \tr es the
 intersection (union) exists and is `\naive'.

\begin{PROP}\label{newt}
Let $(\tT_1,\tT_2)$ be a lower, respectively upper, consistent
pair of \tr es. Then the formulas
$$
\Dlz_\cap=\Dlz_1\cap\Dlz_2,\qquad
\Dgo_\cap=(\Dlz_\cap)^\bot,
$$
with the truncation functor $\tau_{\lle0}^\cap:=\tau_{\lle0}^1\tau_{\lle0}^2$, respectively
$$
\Dgo_\cup=\Dgo_1\cap\Dgo_2,\qquad
\Dlz_\cup={}^\bot(\Dgo_\cup)
$$
with the truncation functor $\tau_{\gge1}^\cup:=\tau_{\gge1}^2\tau_{\gge1}^1$,
define a new \tr e in the category~$\dD$.
\end{PROP}
\begin{PROOF}{.}
For a lower consistent pair $(\tT_1,\tT_2)$ of \tr es the category
$\Dlz_\cap={\Dlz_1\cap\Dlz_2}$ forms the data~1 from lemma~\ref{l1}. Indeed,
it is obviously preserved by the shift functor~$T$. Moreover, the functor
$\tau_{\lle0}^\cap:=\tau_{\lle0}^1\tau_{\lle0}^2$ has its image in $\Dlz_\cap$
due to lower consistency of the pair. One can easily see that it is right
adjoint to the inclusion functor $i_{\lle0}^\cap\colon\Dlz_\cap\to\dD$.
Similarly for an upper consistent pair.
\end{PROOF}

The example from the preceding subsection yields a non-con\-sist\-ent pair
of \tr es.

Now we shall specify conditions on a sequence of \tr es under which one can
iterate the operations of intersection and union.

First, intersection defines an associative operation on ordered
sequences with pairwise lower consistent elements, preserving
lower consistency at intermediate steps.
\begin{PROP}\label{coml}
Let $(\tT_1,\tT_2,\tT_3)$ be a triple of \tr es with $(\tT_i,\tT_j)$
being lower consistent for ${i<j}$. Then the pairs
$(\tT_1\cap\tT_2,\tT_3)$ and $(\tT_1,\tT_2\cap\tT_3)$ are lower consistent.
Moreover, $(\tT_1\cap\tT_2)\cap\tT_3$ coincides with
$\tT_1\cap{(\tT_2\cap\tT_3)}$.
\end{PROP}
\begin{PROOF}{.}
By proposition~\ref{newt} the truncation functor for the \tr e
${\tT_1\cap\tT_2}$ is the composite of $\tau_{\lle0}^1$ and~$\tau_{\lle0}^2$.
As both functors preserve $\Dlz_3$, lower consistency of the pair
$(\tT_1\cap\tT_2,\tT_3)$ follows.

Further, $\Dlz_{2\cap3}=\Dlz_2\cap\Dlz_3$. Both categories $\Dlz_2$ and
$\Dlz_3$ are preserved by~$\tau_{\lle0}^1$. Hence,
$\tau_{\lle0}^1\Dlz_{2\cap3}\subset\Dlz_{2\cap3}$, i.\,e. the second pair
$(\tT_1,\tT_2\cap\tT_3)$ is lower consistent.

Finally, for both \tr es $(\tT_1\cap\tT_2)\cap\tT_3$ and
$\tT_1\cap(\tT_2\cap\tT_3)$ the categories $\Dlz$ coincide with
$\Dlz_1\cap{\Dlz_2\cap\Dlz_3}$. This gives coincidence of these \tr es.
\end{PROOF}

By the dual argument, one proves
\begin{PROP}\label{comr}
Let $(\tT_1,\tT_2,\tT_3)$ be a triple of \tr es with $(\tT_i,\tT_j)$
being upper consistent for ${i<j}$. Then the pairs
$(\tT_1\cup\tT_2,\tT_3)$ and $(\tT_1,\tT_2\cup\tT_3)$ are upper consistent.
Moreover, $(\tT_1\cup\tT_2)\cup\tT_3$ and $\tT_1\cup{(\tT_2\cup\tT_3)}$
coincide.
\end{PROP}

It follows that for any finite sequence $\{\tT_i\}_{i\in[1,n]}$ of
\tr es such that for all ${i<j}$ the pairs $(\tT_i,\tT_j)$ are lower
(resp. upper) consistent, one can construct a new \tr e
$\bigcap\limits_{i\in[1,n]}\tT_i$ (resp. $\bigcup\limits_{i\in[1,n]}\tT_i$)
with $\Dlz_\cap=\bigcap\limits_{i\in[1,n]}\Dlz_i$ (resp.
$\Dgo_\cup=\bigcap\limits_{i\in[1,n]}\Dgo_i$) and the truncation functor
$\tau_{\lle0}^\cap=\tau_{\lle0}^1\tau_{\lle0}^2\ldots\tau_{\lle0}^n$
(resp. $\tau_{\lle0}^\cup=\tau_{\gge1}^n\ldots\tau_{\gge1}^1$).

\bigskip
The next two propositions describe conditions under which one can alternate
the intersection and union operations.

\begin{PROP}\label{ui}
Let $(\tT_1,\tT_2,\tT_3)$ be a triple of \tr es, such that the pairs
$(\tT_1,\tT_2)$ and $(\tT_1,\tT_3)$ are upper consistent and the pair
$(\tT_2,\tT_3)$ is lower consistent. Then the pair $(\tT_1,\tT_2\cap\tT_3)$
is upper consistent.
\end{PROP}
\begin{PROOF}{.}
We denote by $\tau_{\lle0}^\cap$ and $\tau_{\gge1}^\cap$ the truncation
functors for the \tr e ${\tT_2\cap\tT_3}$.

We know by proposition~\ref{newt} that
$\tau_{\lle0}^\cap=\tau_{\lle0}^2\tau_{\lle0}^3$. Then, for any ${X\in\dD}$, by
axioms of triangulated categories, the commutative triangle
\begin{center}
\parbox{50mm}{%
\begin{picture}(50,20)
\put(0,1){$\tau_{\lle0}^2\tau_{\lle0}^3 X$}
\put(11,6){\vector(1,1){11}}
\put(23,18){$\tau_{\lle0}^3 X$}
\put(33,17){\vector(1,-1){11}}
\put(45,1){$X$}
\put(17,2){\vector(1,0){26}}
\end{picture}}
\end{center}
can be embedded into the following `octahedron' with all columns and rows
being exact triangles (for convenience we draw it as a square):
\begin{equation}
\begin{array}{c@{}c@{}c@{}c@{}c}
\tau_{\gge1}^2\tau_{\lle0}^3 X &{}\to{}&\tau_{\gge1}^\cap X&{}\to{}&\tau_{\gge1}^3 X \\
\uparrow & &\uparrow & &{\uparrow}\lefteqn{\wr} \\
\tau_{\lle0}^3 X &{}\to{}&X &{}\to{}&\tau_{\gge1}^3 X \\
\uparrow & &\uparrow & &\uparrow \\
\tau_{\lle0}^2\tau_{\lle0}^3 X &{}\simto{} &\tau_{\lle0}^2\tau_{\lle0}^3 X &{}\to{}& 0
\end{array}
\label{oct1}
\end{equation}

%Since for any \tr e the maps $X_-\to X\to X_+$ in the triangle~(\ref{f1})
%are canonical, the maps in~(\ref{oct1}) are also canonical.

We have to show
that $\tau_{\gge1}^\cap{\Dgo_1\subset\Dgo_1}$. Let ${X\in\Dgo_1}$, then,
since the pair $(\tT_1,\tT_3)$ is upper consistent, then
$\tau_{\gge1}^3X\in\Dgo_1$. From the middle row in~(\ref{oct1}) and remark 2 it follows
that $\tau_{\lle0}^3X\in\Dgo_1$. By upper consistency of the pair
$(\tT_1,\tT_2)$ we deduce that $\tau_{\gge1}^2\tau_{\lle0}^3X\in\Dgo_1$. Now
from the upper row of~(\ref{oct1}) it follows that
$\tau_{\gge1}^\cap X\in\Dgo_1$. This proves the proposition.
\end{PROOF}

Thus, under conditions of proposition~\ref{ui} one can construct a new \tr e
$\tT_1\cup{(\tT_2\cap\tT_3)}$.

By duality, one proves
\begin{PROP}\label{iu}
Let $(\tT_1,\tT_2,\tT_3)$ be a triple of \tr es, such that the pairs
$(\tT_1,\tT_3)$ and $(\tT_2,\tT_3)$ are lower consistent and the pair
$(\tT_1,\tT_2)$ is upper consistent. Then the pair $(\tT_1\cup\tT_2,\tT_3)$
is lower consistent.
\end{PROP}

Suppose we have a pair $(\tT_1,\tT_2)$ of \tr es with
${\tT_1\subseteq\tT_2}$. It is easily seen to be lower and upper consistent.
%As the functor
%$\tau_{\lle0}^2$, while restricted to $\Dlz_2$, is isomorphic to the identity
%functor, it preserves~$\Dlz_1$. Thus
Moreover, $(\tT_2,\tT_1)$ is also a lower and upper consistent pair.
%Moreover, as ${\Dgo_1\supseteq\Dgo_2}$, the same argument
%shows that both pairs $(\tT_1,\tT_2)$ and $(\tT_2,\tT_1)$ are upper
%consistent too.
In particular,
%if $\tT_2=\tT_1[k]$ is obtained from $\tT_1$ by applying the
%shift functor, we conclude that t
the pair $(\tT_1,\tT_1[k])$ is consistent for
any ${k\in\zz}$. Note that truncation functors commute:
$$
\tau^1_{\lle0}\tau^2_{\lle0}=\tau^1_{\lle0}=\tau^2_{\lle0}\tau^1_{\lle0},\ \ \tau^1_{\gge1}\tau^2_{\gge1}=\tau^2_{\gge1}=\tau^2_{\gge1}\tau^1_{\gge1}
$$

\vspace{0.4cm}

We conclude this subsection by proving a lemma, used in the sequel,
 on commutation of truncation
functors for an ordered pair of $t$-structures.

\begin{LEM}\label{comord}

For an ordered pair $(\tT_1,\tT_2)$ of \tr es with
${\tT_1\subseteq\tT_2}$ the truncation functors commute:

$$
\tau^2_{\lle0}\tau^1_{\gge1}=\tau^1_{\gge1}\tau^2_{\lle0}
$$

\end{LEM}

\begin{PROOF}{.}
For any object $X\in \dD$ we have an isomorphism:
$\tau^2_{\gge1}\tau^1_{\gge1}X=\tau^2_{\gge1}X$.

Hence, by the octahedron axiom, we have the following diagram with exact rows
and columns:

\begin{equation}
\begin{array}{c@{}c@{}c@{}c@{}c}
\tau_{\lle0}^2\tau_{\gge1}^1 X &{}\to{}&\tau_{\gge1}^1 X&{}\to{}&\tau_{\gge1}^2
\tau_{\gge1}^1 X \\
\uparrow & &\uparrow & &{\uparrow}\lefteqn{\wr} \\
\tau_{\lle0}^2 X &{}\to{}&X &{}\to{}&\tau_{\gge1}^2 X \\
\uparrow & &\uparrow & &\uparrow \\
\tau_{\lle0}^1 X &{}\simto{} &\tau_{\lle0}^1 X &{}\to{}& 0
\end{array}
\label{oct5}
\end{equation}

In the upper row, two terms  $\tau^1_{\gge1} X$ and $\tau^2_{\gge1}
\tau^1_{\gge1} X$ are in $\Dgo _1$. Hence, the third term $\tau^2_{\lle0}
\tau^1_{\gge1} X$ is in $\Dgo _1$.

Then, in the left column, the lower term is in $\Dlz _1$ and the
upper one is in $\Dgo _1$. Therefore, this column is the canonical
decomposition of the middle term $\tau ^2_{\lle0} X$ along the
first $t$-structure. Hence, we have a canonical isomorphism:
$\tau^2_{\lle0}\tau^1_{\gge1}X\simeq
\tau^1_{\gge1}\tau^2_{\lle0}X.$

\end{PROOF}
%---------------------------------------------------------------------------
\subsection{The standard postulates for lattices}\label{standardpost}
The lattice theory allows to interpret particularly nice posets, called lattices, as abstract algebras with the union (sum) and
intersection (product) operations that satisfy so-called standard postulates~\cite{Bir}. This interpretation is useful for manipulations with lattices, such as taking quotient.

\label{stpo}
\begin{ITEM}{Definition} A {\it lattice\/}
is a partially ordered set for which union and intersection of any pair
of elements exist.
\end{ITEM}

%Unfortunately, I do not know any `large' enough class of \tr es which
%constitute a lattice. This is an important problem to
%find it.
%It would be interesting to be able to find the intersection and
%the union of given elements in such a lattice in an efficient `constructive'
%way (like for consistent pairs).

Recall the {\em standard postulates} known for the operations on a partially
ordered set. We shall use now the (${}\cdot{}$,$+$)-no\-ta\-tion with the usual priority of $\cdot$ over~$+$. In the following formulae, we suppose that both
sides of equations are well-defined for elements $\tT ,\tT_1,\tT_2,\tT_3$ of a partially ordered set (when this is not {\em a~priori} clear, as say for lattices):
\begin{eqnarray}
\lefteqn{\tT\cdot\tT=\tT,\qquad \tT+\tT=\tT;}&& \label{law1}\\
\lefteqn{\tT_1\cdot\tT_2=\tT_2\cdot\tT_1,\qquad
   \tT_1+\tT_2=\tT_2+\tT_1;}&& \label{com}\\
\lefteqn{(\tT_1\cdot\tT_2)\cdot\tT_3=\tT_1\cdot(\tT_2\cdot\tT_3),\qquad
   (\tT_1+\tT_2)+\tT_3=\tT_1+(\tT_2+\tT_3);}&& \label{ass}\\
\lefteqn{\tT_1\cdot(\tT_1+\tT_2)=\tT_1+\tT_1\cdot\tT_2=\tT_1.}&& \label{amb}
\hphantom{wwwwwwwwwwwwwwwwwwwwwwwwwwwwwwwwwwwwwww}
\end{eqnarray}
\begin{PROP}(see \cite{Bir})
Lattices are identified with abstract algebras with two binary operations
$\ \cdot \ $ and $\ + \ $ satisfying standard postulates (\ref{law1})-(\ref{amb}).
Given such an algebra, the partial order on the set of its elements can be recovered from the rule:
\begin{equation}\label{compar}
x\le y \ \ \ {\rm iff} \ \ \ x\cdot y=x.
\end{equation}
If {\bf 0} and/or {\bf 1} exist, they can be added as nullary operations.
\end{PROP}

An equivalence relation $\sim $ in a lattice $L$ is called {\em congruence}
if $p\sim q$ implies $p+l\sim q+l$ and $p\cdot l \sim q\cdot l$ for
any element
$l\in L$. Given a number of pairs (equivalences) $p_i\sim q_i$, $i\in I$, we
define the congruence generated by them as the minimal congruence containing
all of them. It is the intersection of all congruences that contain $\{p_i\sim q_i\}_{i\in I}$.
Intersection of arbitrary number of congruencies (as subsets in $L\times L$)
is again a congruence.

For a congruence $\sim $ in a lattice $L$, regarded as an abstract
algebra, the quotient algebra $L'=L/\sim$ is well-defined. Being an abstract
algebra of the same type (satisfying the same standard postulates), it is
again a lattice.

Note that, for $x',y'\in L$, $x'\le y'$ in $L'$ iff there exist $x, y\in L$ such that $x\sim
x'$, $y\sim y'$ and $x\le y$ in $L$. Indeed, if $x'\le y'$ in $L'$, then
 $x'\sim x'y'$ by (\ref{compar}). Once summed up with $y'$, this gives:
$x'+y'\sim x'y'+y'=y'$. Hence, we can take required $x$ and $y$ to be:
$x=x'y'$, $y=x'+y'$.

Unfortunately, I do not know any `sufficiently large' class of \tr es in a triangulated category which
constitute a lattice. It would be interesting to explore this problem.

\subsection{Modularity laws}

For any partially ordered set, the distributivity inequalities are verified
(again under assumption that both sides are well defined):
\begin{eqnarray}
(\tT_1+\tT_2)\cdot\tT_3&\ge&\tT_1\cdot\tT_3+\tT_2\cdot\tT_3, \label{dis1}\\
\tT_1+\tT_2\cdot\tT_3&\le&(\tT_1+\tT_2)\cdot(\tT_1+\tT_3).\label{dis2}\end{eqnarray}

The distributivity laws, i.\,e. when (\ref{dis1}) and~(\ref{dis2}) are
replaced by the equations with the same left and right hand sides, are not valid for \tr es.

{\bf Example.} Let $\tT_1$ and $\tT_2$ be the \Tn t t-structures
as in the example in subsection \ref{example}. Define a \Tn t \tr e $\tT_3$ by
$$
\Dlz_3=\langle\oO\rangle
$$
Then $\tT_1+\tT_2={\bf 1}$, because $\Dgo_1\cap \Dgo_2=\langle\oO (-3)\rangle\cap \langle\oO \rangle =0$.
Also $\tT_1\cdot \tT_3=\tT_2\cdot \tT_3={\bf 0}$, because $\Dlz_1\cap \Dlz_3= \Dlz_2\cap \Dlz_3=0$. It follows that
the distributivity law
$$
(\tT_1+\tT_2)\cdot\tT_3 = \tT_1\cdot\tT_3+\tT_2\cdot\tT_3
$$
is violated, as the left hand side is $\tT_3$, while the right hand side is ${\bf 0}$.

If $\tT_1\le \tT_3$, then (\ref{dis1}) and (\ref{dis2}) yield the following
autodual inequality:
\begin{equation}
\label{ha}
(\tT_1+\tT_2)\cdot\tT_3\ge \tT_1+\tT_2\cdot\tT_3.
\end{equation}
When replaced by the equation, it is called {\em modularity law.}

\vspace{0.4cm}

We prove {\it
under appropriate assumptions on consistency\/} several statements,
 which might be viewed as modularity laws for \tr es. The first theorem
 is autodual.

\begin{THM}[modularity law 1] \label{ml}
Suppose we are given a triple $(\tT_1,\tT_2,\tT_3)$ of \tr es such that
${\tT_1\le\tT_3}$, the pair $(\tT_1,\tT_2)$ is upper consistent and
the pair $(\tT_2,\tT_3)$ is lower consistent. Then
\begin{itemize}
\item[{\it i})] the pair $(\tT_1+\tT_2,\tT_3)$ is lower consistent,
\item[{\it ii})] the pair $(\tT_1,\tT_2\cdot\tT_3)$ is upper consistent,
\item[{\it iii})] $(\tT_1+\tT_2)\cdot\tT_3=\tT_1+\tT_2\cdot\tT_3$.
\end{itemize}
\end{THM}
\begin{PROOF}{.}
Since ${\tT_1\le\tT_3}$, the pair $(\tT_1,\tT_3)$ is lower consistent.
Hence, lower consistency of the pair $(\tT_1+\tT_2,\tT_3)$ follows from
proposition~\ref{iu}. Similarly, the pair $(\tT_1,\tT_3)$ is upper
consistent, hence $(\tT_1,\tT_2\cdot\tT_3)$ is upper consistent by
proposition~\ref{ui}.

Thus, both sides of {\it iii)\/} are well-defined.

%To show {\it iii)\/}, first, notice that from~(\ref{dis1}) and
%${\tT_1\le\tT_3}$ the inequality follows:
%\begin{equation}
%(\tT_1+\tT_2)\cdot\tT_3\ge\tT_1+\tT_2\cdot\tT_3.
%\label{ha}
%\end{equation}
We have inequality (\ref{ha}).
The opposite to~(\ref{ha}) inequality is, by definition, equivalent to
the inclusion of subcategories:
\begin{equation}
\Dlz_{(1+2)3}\subseteq \Dlz_{1+2\cdot3}. \nonumber\\
%{\vphantom{\bigr)}}^\bot
\end{equation}
%or, in more details,
%\begin{equation}
%{\vphantom{\bigr)}}^\bot \bigl(\Dgo_1\cap\Dgo_2\bigr)\cap\Dlz_3\subseteq
%{\vphantom{\bigr)}}^\bot
%\bigl(\Dgo_1\cap(\Dlz_2\cap\Dlz_3)^\bot\bigr).
%\nonumber
%\end{equation}

This, in turn, is equivalent to the following statement:
$$
\mbox{Let $X\in\Dlz_{(1+2)3}$ and $Y\in\Dgo_{1+2\cdot3}$, then $\Hom(X,Y)=0$.}
\leqno (*)
$$

Let us prove it. To this end, decompose $X$ into an exact triangle:
\begin{equation}
\tau^1_{\lle0}X\to X\to\tau^1_{\gge1}X.
\label{trx}
\end{equation}

As $\tT_1\subseteq(\tT_1+\tT_2)\tT_3$, the pair $(\tT_1,(\tT_1+\tT_2)\tT_3)$
is lower consistent. Hence $\tau^1_{\lle0}X\in\Dlz_{(1+2)3}$. Then, by
remark~\ref{1rem}, $\tau^1_{\gge1}X\in\Dlz_{(1+2)3}$.

Now since $\tT_1\subseteq\tT_1+\tT_2\tT_3$, then
${\Hom(\Dlz_1,\Dgo_{1+2\cdot3})=0}$. It follows that
${\Hom(\tau^1_{\lle0}X,Y)=0}$. Applying the functor ${\Hom({}\cdot{},Y)}$ to
the triangle~(\ref{trx}), we see that it is sufficient to show that
${\Hom(\tau^1_{\gge1}X,Y)=0}$.

Therefore, we have reduced $(*)$ to the case when
$X\in\Dlz_{(1+2)3}\cap\Dgo_1$, ${Y\in\Dgo_{1+2\cdot3}}$.

Now let us decompose $Y$ with respect to the third \tr e:
\begin{equation}
\tau^3_{\lle0}Y\to Y\to\tau^3_{\gge1}Y.
\label{try}
\end{equation}

Since $\tT_{1+2\cdot3}\subseteq\tT_3$, the pair $(\tT_{1+2\cdot3},\tT_3)$ is
upper consistent. Hence $\tau^3_{\gge1}Y\in\Dgo_{1+2\cdot3}$. Then, again by
remark~\ref{1rem}, $\tau^3_{\lle0}Y\in\Dgo_{1+2\cdot3}$.
Since $\tT_{(1+2)3}\subseteq\tT_3$, then ${\Hom(\Dlz_{(1+2)3},\Dgo_3)=0}$. It
follows that ${\Hom(X,\tau^3_{\gge1}Y)=0}$. Applying ${\Hom(X,\cdot)}$
to~(\ref{try}), we see that it is sufficient to show that
${\Hom(X,\tau^3_{\lle0}Y)=0}$.

Therefore, we have reduced $(*)$ to the statement:
$$
\mbox{if $X\in\Dlz_{(1+2)3}\cap\Dgo_1$
and $Y\in\Dlz_3\cap\Dgo_{1+2\cdot3}$, then $\Hom(X,Y)=0$.}
\leqno ({*}{*})
$$

Let us show that
\begin{itemize}
\item[a)] $\Dlz_{(1+2)3}\cap\Dgo_1\subset\Dlz_2$,
\item[b)] $\Dlz_3\cap\Dgo_{1+2\cdot3}\subset\Dgo_2$.
\end{itemize}

Recall that, by proposition~\ref{newt}, since the pair $(\tT_1,\tT_2)$ is
upper consistent, $\tau_{\gge1}^{1+2}=\tau_{\gge1}^2\tau_{\gge1}^1$.
Therefore, if $X\in\Dgo_1$, then $\tau_{\gge1}^{1+2}X=\tau_{\gge1}^2X$. If,
in addition, $X\in\Dlz_{(1+2)3}$, then
$\tau_{\gge1}^2X={\tau_{\gge1}^{1+2}X=0}$, in view of
$\tT_{(1+2)3}\subseteq\tT_{1+2}$. This proves a).

Dually, by proposition~\ref{newt}, since the pair $(\tT_2,\tT_3)$ is
lower consistent, $\tau_{\lle0}^{2\cdot3}=\tau_{\lle0}^2\tau_{\lle0}^3$.
Hence for $Y\in\Dlz_3$, $\tau_{\lle0}^{2\cdot3}Y=\tau_{\lle0}^2Y$. If,
in addition, $Y\in\Dgo_{1+2\cdot3}$, then
$\tau_{\lle0}^2Y={\tau_{\lle0}^{2\cdot3}Y=0}$, in view of
$\tT_{2\cdot3}\subseteq\tT_{1+2\cdot3}$. This proves b).

Since a) and b) obviously imply~$(**)$, this finishes the proof of the
theorem.
\end{PROOF}

We shall also need another incarnation of the modularity law.

\begin{THM}[modularity law 2] \label{ml2}
Suppose we are given a triple $(\tT_1,\tT_2,\tT_3)$ of \tr es such that
${\tT_1\le\tT_3}$, the pair $(\tT_1,\tT_2)$ is upper consistent and
the pair $(\tT_3,\tT_2)$ is lower consistent. Then
\begin{itemize}
\item[{\it i})] the pair $(\tT_1,\tT_2\cdot\tT_3)$ is upper consistent,
\item[{\it ii})] the pair $(\tT_3, \tT_1+\tT_2)$ is lower consistent,
\item[{\it iii})] $(\tT_1+\tT_2)\cdot\tT_3=\tT_1+\tT_2\cdot\tT_3$.
\end{itemize}
\end{THM}

Note that in the contrary to the theorem \ref{ml}, lower consistency of
the pair
 $(\tT_3, \tT_1+\tT_2)$ is not a formal consequence of the propositions from
the previous subsection.

\vspace{0.4cm}

\begin{PROOF}{.}
As $\tT_1\le \tT_3$, the pair $(\tT_1,\tT_3)$ is upper consistent. Hence,
we can apply proposition \ref{ui} to the triple $(\tT_1, \tT_3,\tT_2)$
to get upper consistency of the pair $(\tT_1, \tT_2\cdot\tT_3)$.
This proves {\it i)} and the fact that $\tT_{1+2\cdot 3}$ exists.

We know by proposition \ref{newt} that $\tau^{2\cdot 3}_{\lle0}=
\tau^3_{\lle0} \tau^2_{\lle0}$.

Then, for any $X\in \dD$, by the octahedron axiom, we have the
diagram with exact rows and columns:

\begin{equation}
\begin{array}{c@{}c@{}c@{}c@{}c}
\tau_{\gge1}^3\tau_{\lle0}^2 X &{}\to{}&\tau_{\gge1}^{2\cdot 3} X&{}\to{}&\tau_{\gge1}^2 X \\

\uparrow & &\uparrow & &{\uparrow}\lefteqn{\wr} \\

\tau_{\lle0}^2 X &{}\to{}&X &{}\to{}&\tau_{\gge1}^2 X \\

\uparrow & &\uparrow & &\uparrow \\

\tau_{\lle0}^3 \tau_{\lle0}^2 X &{}\simto{} &\tau_{\lle0}^3 \tau_{\lle0}^2 X &{}\to{}& 0
\end{array}
\label{oct2}
\end{equation}
By proposition \ref{newt}, $\tau_{\gge1}^{1+2}=\tau_{\gge1}^2\tau_{\gge1}^1$
and $\tau_{\gge1}^{1+2\cdot 3}=\tau^{2\cdot 3}_{\gge1}\tau_{\gge1}^1$.
Hence, substitution $X=\tau_{\gge1}^1Y$ in the upper row of the octahedron
yields:
$$
\tau^3_{\gge1}\tau^2_{\lle0}\tau^1_{\gge1}Y\to \tau_{\gge1}^{2\cdot 3+1}Y\to
\tau^{1+2}_{\gge1}Y.
$$

Now, substitution $Y=\tau_{\lle0}^{1+2}Z$ gives a natural isomorphism:
$$
\tau^3_{\gge1}\tau^2_{\lle0}\tau^1_{\gge1}\tau^{1+2}_{\lle0}Z\simto
\tau_{\gge1}^{1+2\cdot 3}\tau_{\lle0}^{1+2}Z.
$$

It follows that, for any $Z\in \dD$,
$$
\tau^3_{\lle0}\tau^{1+2\cdot 3}_{\gge1}\tau^{1+2}_{\lle0}Z=0.
$$

Observe an obvious inequality: $\tT_{1+2\cdot 3}\le \tT_3$. Hence, by lemma
\ref{comord}, $\tau_{\lle0}^3$ and $\tau _{\gge1}^{1+2\cdot 3}$ commute. This
implies:
$$
\tau_{\gge1}^{1+2\cdot 3}\tau_{\lle0}^3\tau_{\lle0}^{1+2}=0.
$$

Equivalently, the image of $\tau^3_{\lle0}\tau^{1+2}_{\lle0}$ is in $\Dlz
_{1+2\cdot 3}$.

On the other hand, $\Dlz_{1+2\cdot 3}\subset \Dlz_{1+2}\cap \Dlz_3$. Therefore,
the image of $\tau^3_{\lle0}\tau^{1+2}_{\lle0}$ is in $\Dlz_{1+2}\cap \Dlz_3$.
Then $\tau _{\lle0}^3\Dlz_{1+2}\subset \Dlz_{1+2}$, hence lower consistency of
$(\tT_3,\tT_1+\tT_2)$.

By proposition \ref{newt}, functor $\tau^3_{\lle0}\tau^{1+2}_{\lle0}$ is the
truncation functor for $\tT_{(1+2)\cdot 3}$. Its image is $\Dlz _{(1+2)\cdot
3}=\Dlz_{1+2}\cap \Dlz_3$. It was shown that the image is in $\Dlz_{1+2\cdot
3}$. Together with the opposite inclusion (\ref{ha}) this yields {\it iii)}.
\end{PROOF}

The dual statement to theorem \ref{ml2} is

\begin{THM}[modularity law 2${}^{\prime}$] \label{ml2p}
Suppose we are given a triple $(\tT_1,\tT_2,\tT_3)$ of \tr es such that
${\tT_1\le\tT_3}$, pair $(\tT_2,\tT_1)$ is upper consistent and
pair $(\tT_2,\tT_3)$ is lower consistent. Then
\begin{itemize}
\item[{\it i})]  pair $(\tT_2\cdot\tT_3,\tT_1)$ is upper consistent,
\item[{\it ii})]  pair $(\tT_1+\tT_2, \tT_3)$ is lower consistent,
\item[{\it iii})] $(\tT_1+\tT_2)\cdot\tT_3=\tT_1+\tT_2\cdot\tT_3$.
\end{itemize}
\end{THM}

\section{Sets with consistencies and pairs of chains}

For the sake of clear and logical exposition, we formalize the properties
of consistency proved in the previous subsections and present a convenient
pictorial description.
Then we prove that a consistent pair of chains generates a distributive lattice.

\subsection{ Sets with consistencies}
\label{setc}

By a {\em set with consistencies} we call a partially ordered set containing ${\bf
0}$ and ${\bf 1}$ with two binary
(non-reflexive) relations, called lower and upper consistencies, such that
for any lower, respectively upper, consistent pair the abstract intersection,
respectively union, exists and satisfies the following axioms ( {\em
prime} means the dual axiom):
\begin{itemize}
\item[$SC1.$] if $a\le b$, then both $(a,b)$ and $(b,a)$ are lower and upper
consistent,
\item[$SC2.$] implication of lower consistency under intersection as in proposition \ref{coml},
\item[$SC2^{\prime}.$] implication of upper consistency under union as in
proposition \ref{comr},
\item[$SC3.$] implication of upper consistency under intersection as in
proposition \ref{ui},
\item[$SC3^{\prime}.$] implication of lower consistency under union as in
proposition
\ref{iu},
\item[$SC4.$] modularity law as in theorem \ref{ml},
\item[$SC5.$] modularity law as in theorem \ref{ml2},
\item[$SC5^{\prime}.$] modularity law as in theorem \ref{ml2p}.

\end{itemize}

In this formalism, the words `of $t$-structures' must be removed from all
the quoted statements. In the sequel, we shall sometimes directly refer to the
corresponding statements, meaning the above axioms for general sets with
consistencies.

Note that the standard postulates (\ref{law1}),(\ref{com}),(\ref{ass}),(\ref{amb})
are automatically satisfied
for sets with consistencies. It follows from $SC1$ that {\bf 0} and {\bf 1} constitute a lower and upper
consistent pair with any element.

Contrary to a lattice, a set with consistencies is not an abstract algebra.

%\vspace{0.5cm}

\subsection{Diagrams}
\label{pict}
It is instructive to use
the following pictorial description for consistencies
and operations with them.

Suppose we are given a set with consistencies $V$. Consider a graph $\Gamma
(V)$ with vertices labelled by elements of $V$.

We draw a full arrow:
$$
a\longrightarrow b
$$
if $(a,b)$ is a lower consistent pair in $V$, and a dashed arrow
$$
a \dashrightarrow b
$$
if $(a,b)$ is an upper consistent pair.

Now the operations of union and intersection are interpreted
 as contractions
of full and dashed arrows. For the intersection, this is:
$$
a\longrightarrow b \ \ \ \Longrightarrow \ \ \ a\cdot b
$$

For the union, this is:
$$
a\dasharrow b \ \ \ \Longrightarrow \ \ \ a+b
$$

The duality on $t$-structures described in subsection \ref{poset} exchanges
full arrows with dashed ones and reverses their directions.

We shall actually consider subgraphs of $\Gamma (V)$ which are significant to ensure
iterated contractions.

Axioms $SC2$, $SC2'$, $SC3$, $SC3'$ of the definition
of a set with consistencies provides us with patterns which allow to connect
by arrows a new vertex obtained as result of contractions of an arrow with
other (old) vertices of the graph. Axiom $SC2$ (proposition \ref{coml})
reads as possibility of contractions of two edges in a triangle with full
arrows:

\begin{center}
$a\longrightarrow b\cdot c\ \ \  \Longleftarrow$\ \ \ \ \
\parbox{40mm}{%
\begin{picture}(35,20)
\put(0,1){$a$}
\put(4,4){\vector(1,1){11}}
\put(16,16){$b$}
\put(19,15){\vector(1,-1){11}}
\put(31,1){$c$}
\put(4,2){\vector(1,0){25}}
\end{picture}}
$\Longrightarrow \ \ \ a\cdot b \longrightarrow c$
\end{center}

Axiom $SC2'$ reads as the same picture with dashed arrows.

Axiom $SC3$ (proposition \ref{ui}) reads:

\begin{center}
\parbox{40mm}{%
\begin{picture}(35,20)
\put(0,1){$a$}
\put(4,4){\line(1,1){4}}
%\put(8,8){\line(1,1){3}}
\put(10,10){\vector(1,1){5}}
\put(16,16){$b$}
\put(19,15){\vector(1,-1){11}}
\put(31,1){$c$}
\put(4,2){\line(1,0){6}}
\put(12,2){\line(1,0){6}}
\put(20,2){\vector(1,0){9}}
\end{picture}}
$\Longrightarrow \ \ \ a \dashrightarrow b\cdot c$
\end{center}

Axiom $SC3'$ reads:

\begin{center}
\parbox{40mm}{%
\begin{picture}(35,20)
\put(0,1){$a$}
\put(4,4){\line(1,1){4}}
%\put(8,8){\line(1,1){3}}
\put(10,10){\vector(1,1){5}}
\put(16,16){$b$}
\put(19,15){\vector(1,-1){11}}

\put(31,1){$c$}
\put(4,2){\vector(1,0){25}}
\end{picture}}
$\Longrightarrow \ \ \ a+b \longrightarrow c$
\end{center}

Note that we don't have in general the implications which are
obtained from the last two pictures by changing full arrows for dashed
ones and vice versa (while preserving their directions).

Sometimes it is useful to put the vertex marked with the intersection
(resp. union) \tr e at  the middle point of the corresponding full (resp. dashed) arrow,
if it exists.
%\vspace{0.5 cm}

\subsection{The universal lattice with consistencies}
\label{universal}
A morphism of sets with consistencies is a morphism of the corresponding
posets which takes any upper (lower) consistent pair to an upper (lower) consistent
pair.
A set with consistencies which is a lattice is called a {\em lattice
with consistencies}. Morphisms of lattices with consistencies are defined
in the obvious way.

A morphism $\phi :P\to L$, where $P$ is a poset and $L$ a lattice, is an order
preserving map such that $\phi (x+y)=\phi (x) + \phi (y)$, $\phi (xy)=\phi
(x)\phi (y)$, whenever $x+y$ or $xy$ exists.

Consider a poset $P$ on which two arbitrary binary relations are fixed. We
are interested in morphisms $P\to L$ of $P$ into a lattice $L$ with consistencies,
such that the images of any pair of elements in $P$ which is in the first,
respectively, the second, relation, is lower, respectively upper, consistent
in $L$. Such a morphism $P\to U$ is called {\em universal} if for any other
morphism $P\to L$ of this sort, there exists a unique morphism of lattices
$U\to L$ with consistencies which makes commutative the following diagram
(of morphisms):

\begin{center}
\parbox{40mm}{%
\begin{picture}(35,20)
\put(0,1){$P$}
\put(4,4){\vector(1,1){11}}
\put(16,16){$U$}
\put(19,15){\vector(1,-1){11}}
\put(31,1){$L$}
\put(4,2){\vector(1,0){25}}
\end{picture}}
\end{center}

In this case, we call $U=U(P)$ the {\em universal lattice with consistencies}
generated by $P$.

When consistency relations are omitted in the above definitions,
the notion of the {\em universal lattice}
$L(P)$ generated by a poset $P$ takes place. The universal lattice $L(P)$ classically known to
exist \cite{Bir}. Indeed, remember that the notion of a lattice is equivalent
to that of abstract algebra satisfying the standard postulates. Thus, we can
consider a free abstract algebra $F(P)$, generated by $P$ as a disjoint set,
and take its quotient by the congruence generated by relations $x+y="x+y"$ and $x\cdot y= "x\cdot y"$
for all $x$ and $y$ for which the union $"x+y"$ (resp. intersection $"x\cdot
y"$) in  poset $P$ exists (cf. \ref{standardpost}). This quotient is $L(P)$.

A lattice is said to be
{\em distributive}
if for any $x,y,z$:
$$
x(y+z)=xy+xz,
$$
$$
x+yz=(x+y)(x+z).
$$
These postulates are mutually dual and each one implies the other.

We can similarly define the {\em universal distributive lattice} $D(P)$ generated
by poset $P$. $D(P)$ is the quotient of $L(P)$ by the congruence generated by
equivalences $(x+y)z\sim xz+yz$ for all $x,y,z\in L(P)$.

\begin{PROP} \label{unilat}
Let $P$ be a partially ordered set with two binary relations. Then a universal lattice with consistencies generated by $P$ exists.
\end{PROP}

\begin{PROOF}{.} Consider the universal lattice $L(P)$ generated by the poset
$P$.
Endow $L(P)$ with two binary relations, called respectively upper and lower
consistency: first, demand the images in $L(P)$ of all pairs in $P$ which
belong to the first (resp. second) given binary relation to be upper (resp.
lower) consistent; second, extend the two binary relations of consistency by the minimal set of pairs
that satisfy axioms $SC1$,\dots , $SC3'$ and those parts of axioms $SC5$ and $SC5'$ which
imply consistencies (recall that not all of them formally follow from
the other axioms).

Note that this procedure does not make $L(P)$ a set
with consistencies, because the modularity equations of axioms $SC4$, $SC5$,
$SC5'$ are not satisfied yet.

To meet this constraints, consider the quotient lattice $L_1(P)$ of $L(P)$
by the congruence generated by equivalences $(t_1+t_2)t_3\sim t_1+t_2t_3$, for
all the triples $(t_1,t_2,t_3)\in L(P)$ which meet conditions of one of the modularity
laws $SC4$, $SC5$, $SC5'$.
Endow $L_1(P)$ with consistency relations inherited from $L(P)$. The axiom
$SC1$ obviously follows from the fact (proved in \ref{stpo}) that $x'\le y'$ in
$L_1(P)$ iff there exists $x,y \in L(P)$ such that $x\sim x'$, $y\sim y'$
and $x\le y$ in $L(P)$.

But it is not so good in what concerns the other axioms. For this reason, subject
$L_1(P)$ (with its consistency relations) to the same procedure as
we did with $L(P)$, i.e. extend appropriately the binary relations and take the similar quotient,
to get $L_2(P)$ from $L_1(P)$. By iterating this process, we get a sequence of lattices $L_i(P)$
together with maps $P\to L_i(P)$, compatible with the quotient maps.
Then $U(P)={\rm colim}\ L_i(P)$  with two colimit binary relations is the universal lattice with consistencies.

Indeed, consider any morphism $P\to L$, where $L$ is a lattice with consistencies
and the morphism takes the pairs in $L$ from the first, respectively, the second binary relations on $P$ into lower, respectively, upper consistent pairs in $L$. The universality of $P\to L(P)$ implies existence of the lattice homomorphism $L(P)\to L$, which makes the corresponding diagram commute. This homomorphism clearly descends to a homomorphism $L_i(P)\to L$, which, when passing to the limit, gives a homomorphism $U(P)\to L$.
It is easy to see that this homomorphism is unique.

\end{PROOF}

We can consider a distributive lattice as a lattice with consistencies for which all pairs of elements are both upper and lower consistent. Since distributivity implies modularity, we indeed get a correctly defined lattice with consistencies.
A canonical homomorphism
$$
\psi :U(P)\to D(P)
$$
follows from the universality of $U(P)$. It is clearly an epimorphism, because congruencies used to construct $U(P)$ as a quotient of the lattice $L(P)$ in proposition \ref{unilat} are among those which were used to construct $D(P)$ as the quotient of $L(P)$.

\subsection{Consistent pairs of chains}

A {\em chain} in a partially ordered set is an ordered sequence of elements.
A pair of chains $a_1\ge \dots \ge a_n$ and
$b_1\le \dots \le b_m$ in a set with consistencies
such that $(a_i, b_j)$ is upper and lower consistent for any pair $(i,j)$,
is called {\em consistent pair of chains}. Our purpose in this subsection
will be to prove that a consistent pair of chains generates (under operations
of union and intersection) a
distributive lattice.

There is a theorem due to Birkhoff \cite{Bir} that two chains in a modular
lattice generate a distributive lattice. When working with sets with consistencies we have two new difficulties
against the hypothesis of Birkhoff theorem: first, intersection and union
do not {\em a priori} exist, second, the modularity laws are verified only under some
consistency conditions.

\vspace{0.4cm}

The notion of sum and intersection of any number of elements in a partially
ordered set is intrinsically defined. Say, the sum of a given set $S$ of elements
is the element which is greater than any element from $S$ but less than
any other element with this property. It may not exist.
Recall the following trivial observation, which slightly strengthens
associativity postulates (\ref{ass}).

\begin{LEM}\label{strass}
Let $(a,b,c)$ be a triple of elements in a partially ordered set. If $a+b$
exists, then existence of one side of the following equation implies existence
of the other and the equation itself:

\begin{equation}\label{abcabc}
(a+b)+c=a+b+c.
\end{equation}
In the dual statement, + is replaced by $\cdot $.
\end{LEM}

\begin{PROOF}{.} Obvious.
\end{PROOF}
\vspace{0.4cm}

{\bf Remark.} Note that in the constructions of this subsection, this lemma is the only
tool applied to deduce existence of union, respectively intersection,
 for pairs of elements
($(a+b, c)$ here) which do not {\em a priori} constitute an upper,
respectively lower, consistent pair. The typical situation is the following.
Suppose we have a triple of elements (say $t$-structures) $(a_1,a_2,a_3)$
with $a_i\to a_j$, for $i<j$. Then, by proposition \ref{iu} and the lemma,
$(a_1+a_2)+a_3=a_1+a_2+a_3$ exists. Hence, again by the lemma, $(a_1+a_3)+a_2$
exists, though neither $(a_1+a_3,a_2)$ nor $(a_2,a_1+a_3)$ is {\em a priori}
lower consistent.

\vspace{0.4cm}

For the rest of the subsection, we assume to be given a consistent pair of
chains $a_1\ge \dots \ge a_n$ and
$b_1\le \dots \le b_m$.
Denote:
$$
u_{ij}=a_ib_j,\ \ \  v_{ij}=a_i+b_j.
$$

\begin{LEM}\label{upu}
Let $j\le l$. Then, the pair $(u_{kl},u_{ij})$ is upper consistent.
Dually, the pair $(v_{ij},v_{kl})$ is lower consistent for $i\le k$.
\end{LEM}

\begin{PROOF}{.} First, we may suppose that
$i\le k$, because, if $i>k$, then $u_{ij}\le u_{kl}$ and upper consistency
for $(u_{kl},u_{ij})$ follows.

Applying $SC3$
%proposition \ref{ui}
to $(a_k,a_i,b_j)$, we obtain upper consistency
for the pair $(a_k, a_ib_j)$.

Hence,
$SC5'$ (the modularity law $2'$) is applicable to the triple $(a_ib_j,a_k,b_l)$.
It follows, the pair $(u_{kl},u_{ij})=(a_kb_l,a_ib_j)$ is upper consistent.
The rest follows from the duality.
\end{PROOF}

\vspace{0.5cm}

Recall that
a decomposition $x=x_1+\dots +x_n$ is called {\em irreducible}, if $x\ne
x_1+\dots
+x_{i-1}+x_{i+1}+\dots +x_n$, for all $i$.

\begin{PROP}\label{sumu}
Any finite sum of elements $u_{ij}$ exists and has a decomposition of the
form:
\begin{equation}\label{uijk}
u=u_{i_1j_1}+\dots +u_{i_kj_k}, \ with \ i_1<\dots <i_k; \ j_1<\dots <j_k.
\end{equation}

Dually, any finite product of $v_{ij}$ exists and has the form
\begin{equation}\label{vijk}
v=v_{i_1j_1}\cdot \dots \cdot v_{i_kj_k}, \ with \ i_1<\dots <i_k; \ j_1<\dots
<j_k.
\end{equation}
\end{PROP}
\begin{PROOF}{.}
First, any sum of the form (\ref{uijk}) exists. Indeed, summation can be
done step by step from the left to the right in view of lemma \ref{upu} and
$SC2'$.
%proposition \ref{comr}.

The rest of the proof is done by induction on number of summands $u_{ij}$.

Let us add to a sum $u$ of the form (\ref{uijk}) some element $u_{ij}$. In
fact, it is sufficient to prove only existence of
\begin{equation}\label{uprim}
u'=u+u_{ij}.
\end{equation}
Indeed, if it
exists, then we can reduce, if necessary, the decomposition of $u'$ to get
an irreducible one. Any irreducible decomposition is automatically of the
form (\ref{uijk}), because otherwise there would exist summands $u_{pq}$
and $u_{rs}$ with $u_{pq}\le u_{rs}$.

If $i\le i_t$ and $j\ge j_t$ for some $t$, then $u_{ij}\ge u_{i_tj_t}$.
Similarly, if $i\ge i_t$ and $j\le j_t$ for some $t$, then $u_{ij}\le
u_{i_tj_t}$. In both cases, one of two elements in the sum for $u'$, either
$u_{ij}$ or $u_{i_tj_t}$ can be absorbed by the rule (\ref{amb}). Then, in
view of lemma \ref{strass} and the induction hypothesis the sum $u'$ exists.

The alternative to the above inequalities is the possibility to order all
summands in the decomposition (\ref{uprim}) of $u'$ to be of the form
(\ref{uijk}). It has already been shown that any such sum exists.

\end{PROOF}

A chain is called {\em extended} if it contains ${\bf 0}$
and ${\bf 1}$ as extreme elements.
\begin{PROP} \label{mulesum}
Let $a_1\ge \dots \ge a_n$ and $b_1\le \dots \le b_m$ be
a consistent pair of extended chains. Choose two ordered sets of indices $I=\{i_1\le \dots \le i_k\}$, $i_l\in [1,n]$, and $J=\{j_1 \le \dots \le j_k\}$, $j_l\in [1,m]$. Denote:
$$
r_{IJ}=a_{i_1}(b_{j_1}+a_{i_2})\dots(b_{j_{k-1}}+a_{i_k})b_{j_k},
$$
$$
s_{IJ}=a_{i_1}b_{j_1}+a_{i_2}b_{j_2}+\dots +a_{i_k}b_{j_k}.
$$
Then:
\begin{itemize}
\item[{\em i)}] $r_{IJ}=s_{IJ}$,
\item[{\em ii)}] $a_i\dasharrow s_{IJ}\longrightarrow b_j$, for all $(i,j)\in
[1,n]\times[1,m]$.
\end{itemize}
\end{PROP}
{\bf Remark.} If we extend the chains of $\{a_i\}_{i\in I}$ and $\{b_j\}_{j\in J}$ by ${\bf 0}$, then the formula for $r_{IJ}$ is simply the product of suitable $v_{ij}$'s. We allow the indices for $\{a_i\}$'s and $\{b_j\}$'s coincide for the reason of simplifying the write-down of the proof by induction.

\begin{PROOF}{.} Both $r_{IJ}$ and $s_{IJ}$ exist due to
proposition \ref{sumu}.

We shall prove both statements of the proposition by
simultaneous induction on $k=|I|=|J|$. For $k=1$,
$r_{IJ}=a_{i_1}b_{j_1}=s_{IJ}$
and {\em ii)} follows from $SC2$ and $SC3'$.
%propositions \ref{ui} and \ref{coml}.

Now, take the case $|I|=|J|=k-1$ for granted.

Denote $I'=I\setminus i_k$, $J'=J\setminus j_k$. Evidently, $s_{I'J'}\le b_{j_k}$.
By the induction
hypothesis $a_{ij}\dasharrow s_{I'J'}$. Then, we can apply $SC5'$
%modularity law $2'$
(theorem \ref{ml2p}) to the triple $(s_{I'J'},a_{i_k},b_{j_k})$ to
get
\begin{equation}\label{sabr}
(r_{I'J'}+a_{i_k})b_{j_k}=(s_{I'J'}+a_{i_k})b_{j_k}=s_{I'J'}+a_{i_k}b_{i_k}=s_{I
J}.
\end{equation}

Let us show that the left hand side is $r_{IJ}$.

Decompose $r_{I'J'}=qb_{j_{k-1}}$, where
$$
q=a_{i_1}(b_{j_1}+a_{i_2})\dots (b_{j_{k-2}}+a_{i_{k-1}}).
$$
Since the second chain is extended, the last element $b_m$ in it is ${\bf 1}$.
Then, $q=s_{I'J''}$, where
$J''=\{j_1,\dots j_{k-2},m\}$. Hence, applying the induction hypothesis
to the subsets $I'$
and $J''$, we have: $q\longrightarrow b_{j_{k-1}}$.
Since, in addition,  $a_{i_k}\le q$ and $a_{i_k}\dasharrow b_{j_{k-1}}$,
we can apply $SC5$
%modularity law 2
to the triple $(a_{i_k},b_{j_{k-1}},q)$:
$$
r_{I'J'}+a_{i_k}=qb_{j_{k-1}}+a_{i_k}=q(b_{j_{k-1}}+a_{i_k}).
$$

It follows that the left hand side of (\ref{sabr}) is $r_{IJ}$. Thus, we
have proven {\em i)}.

In view of $a_{i_k}\dasharrow s_{I'J'}$ and the induction hypothesis, we
can apply axioms $SC2'$ and $SC3'$
%\ref{comr} and \ref{iu}
to get:
$$
a_i\dasharrow \ (s_{I'J'}+a_{i_k})\longrightarrow b_j
$$
for all $(i,j)$. Then, in view of $s_{I'J'}+a_{i_k}\longrightarrow b_{j_k}$
we can apply axioms
$SC2$ and $SC3$
%\ref{ui} and \ref{coml}
to get:
$$
a_i\dasharrow (s_{I'J'}+a_{i_k})b_{j_k}\longrightarrow b_j.
$$
By (\ref{sabr}) the middle entry coincides with $s_{IJ}$.
This proves {\em ii)} and completes the proof of the proposition.
\end{PROOF}

\begin{THM}\label{latti}
A consistent pair of chains in a set with consistencies generates in this set a distributive
lattice with all elements
being of the form (\ref{uijk}).
\end{THM}
\begin{PROOF}{.}
%An arbitrary pair of chains also generates a lattice,
%because we can extend them
By adding, if necessary, ${\bf 0}$ and ${\bf 1}$ as extreme elements of the chains,
we may assume
%For the beginning, suppose that
chains to be extended.
Then all $a_i$ and $b_j$ are among $u_{ij}$.  By
proposition \ref{sumu}, any sum of $u_{ij}$ exists and is of the form
(\ref{uijk}).

Let us show that any product of elements of the form (\ref{uijk})
exists and is of the same form. By proposition \ref{mulesum}, any sum of the
form (\ref{uijk}) can be recast as a product of $v_{ij}$ (in the formula
for $r_{IJ}$ elements $a_{i_1}$ and $b_{j_k}$ can be interpreted as $v_{i_1m}$
and $v_{1j_k}$ respectively).
Products of these elements exist again by proposition \ref{sumu} and are of
the form (\ref{vijk}), which in turn, by proposition \ref{mulesum}, can
be recast as a sum of the form (\ref{uijk}). This proves that the set of
elements of the form (\ref{uijk}) is a lattice, $L$. It is the lattice generated
by the pair of chains.
%An arbitrary pair of chains also generates a lattice,
%because we can extend them
%by adding ${\bf 0}$ and ${\bf 1}$ as extreme elements of the chains.

Consider the poset $P=\sigma \cup \tau$, which is a disjoint union of two
chains $\sigma =\{a_1\ge \dots \ge a_n\}$ and
$\tau =\{b_1\le \dots \le b_m\}$, as an abstract ordered set. Let $D(P)$ be
the universal distributive lattice generated by $P$.
%The quotient of it by the distributivity equations is the
%{\em universal distributive lattice} $D(P)$ generated by $P$.

Endow  $P$ with two coinciding binary relations, which consist of all pairs
$(a_i, b_j)$. Let $U(P)$ be the universal lattice with consistencies generated
by $P$, which exists by proposition \ref{unilat}, and $\psi :U(P)\to D(P)$ the
canonical epimorphism. Let us show that it is an isomorphism.

Since $U(P)$ is the quotient of the universal lattice $L(P)$, it has only elements of the form (\ref{uijk}).
Hence, to show that $\psi $ is an isomorphism it is sufficient to present a
distributive lattice generated by two chains such that all elements of the
form (\ref{uijk}) are taken by $\psi $ into distinct elements.

Such a lattice can be realized as a sublattice of subsets in the
finite set of integer nodes in the rectangle  $\{(x,y)\in
[1,m]\times [-1,-n]\}$. Elements $a_i$ are presented by the
subsets $\{y\le -i\}$ and $b_j$ by the subsets $\{x\le j\}$. Any
element of the lattice is the set of integer points under a ladder
descending from the point $(1,-1)$ to the point $(m,-n)$. Being a
sublattice of the lattice of subsets, it is distributive and
images under $\psi $ of distinct elements of the form (\ref{uijk})
are clearly distinct. Hence $\psi $ is injective. Therefore,
$U(P)$ is distributive and, by universality property, $L$ is
distributive too.

\end{PROOF}

\vspace{0.4cm}

{\bf Remark.}
An element $x$ in a lattice
is called {\em indecomposable} if $x=a+b$ implies $x=a$ or $x=b$, and {\em
decomposable} otherwise.
By another theorem of Birkhoff \cite{Bir}, any element in a finite
distributive lattice can be uniquely decomposed into an irreducible sum of
indecomposable elements. Clearly, all indecomposable elements are among
$u_{ij}$. It is easy to see that an element $u_{ij}$ is decomposable in $L$
iff $u_{ij}=u_{i+1,j}+u_{i,j-1}$.

\vspace{0.4cm}

\section{Perverse coherent sheaves}
\subsection{Grothendieck duality and the dual \tr e}
Let $k$ be a field of characteristic zero. In this section we
assume, for simplicity, all schemes to be of finite type over $k$.
For such a scheme $X$ we consider the bounded derived category
$\dD (X):=\dD^b_{\rm coh}(X)$ of coherent sheaves. We shall
construct some perverse \tr es in $\dD (X)$ using the machinery
developed in the preceding section.

Let $\tT $
%= (\Dlz ,\Dgo )=(\Dlz (X), \Dgo (X))$
be the standard \tr e in
$\dD (X)$ with $\Dlz $ (resp. $\dD ^{\ge 0}$) consisting of complexes of coherent
sheaves on $X$ with trivial positive (resp. negative) degree cohomology.

By ${\cal H}^i ({\cal F})$ we denote the $i$-th cohomology sheaf of a complex
${\cal F}\in \dD (X)$.
${\cal O}_X$ is the
structure sheaf of rings on $X$.
We write $f_*$, ${\cal H}om$, etc.
for derived push-forward, derived local Hom functor and other functors between derived categories. For a functor
$\Phi :{\cal D}(X)\to {\cal D}(Y)$ we use notation ${\mathbb R}^k\Phi :={\cal H}^k\cdot \Phi $
for the corresponding cohomological functors ${\cal D}(X)\to {\rm Coh}(Y)$. Here $\Phi $ is usually the derived functor
of a left exact functor ${\rm Coh}(X)\to {\rm Coh}(Y)$. If $\Phi $ is the derived one of a right exact functor
${\rm Coh}(X)\to {\rm Coh}(Y)$,
%(like $i^*$ for a closed embedding $i:Y\to X$),
then, in compliance with the tradition, we write $\mathbb{L}^k\Phi :={\cal H}^{-k}\cdot \Phi $
for the cohomology functors.

Note that for a morphism $f:X\to Y$ there is the twisted inverse image functor
$f^!:\dD _{qcoh}^+(Y)\to \dD_{qcoh}^+(X)$ between the bounded below derived categories
of quasi-coherent sheaves (see \cite{V1} and \cite{LM}). It is defined by:
$$
f^!(-)=(-)\otimes f^!({\cal O}_Y),
$$
where $f^!({\cal O}_Y)$ is defined by applying the right adjoint functor for $f_*$ to the structure sheaf ${\cal O}_Y$.
When $f$ is proper and of finite Tor-dimension, then $f^!$ coincides with the
right adjoint to $f_*$ and $f^!$ takes $\dD (Y)$
to $\dD (X)$. In the particular case when $f$ is a closed embedding, this property gives an obvious construction for $f^!$ at least locally over the base. When $f$ is smooth of relative dimension $n$, $f^!$ is defined by twisting with the relative canonical class $\omega _{X/Y}$:
$$
f^!(-)=(-)\otimes \omega _{X/Y}[n].
$$
If $f$ is projective, it can be decomposed into $f=pi$, where $i$
is a closed embedding and $p$ is smooth. Then $f^!=p^!i^!$. Note also that $f$ is automatically of finite Tor-dimension if $Y$ is a smooth variety.

Also, for a proper $f$, there is a
duality isomorphism for local ${\cal H}om$'s:
\begin{equation}\label{adish}
f_*{\cal H}om_X({\cal F},f^!{\cal G})={\cal H}om_Y(f_*{\cal F},{\cal G}),
\end{equation}
natural in  ${\cal F}\in {\cal D}^+_{qcoh}(X)$, ${\cal G}\in {\cal D}^+_{qcoh}(Y)$.
%In particular, it works for
%${\cal F}\in {\cal D}(X)$, ${\cal G}\in {\cal D}(Y)$.

Denote by $\omega ^\bcdot _X$ the dualizing complex on $X$.
% which is known to be an object in $\dD^b_{\rm coh}(X)$
%for $X$ of finite type over $k$
 By definition $\omega ^\bcdot _X=\pi ^!_X{\cal O}_{\rm pt}$, where $\pi_X:X\to {\rm pt}$
is  the projection to a point pt. For a smooth variety $X$ of dimension
$n$, $\omega ^{\bcdot }_X=\omega _X[n]$, where $\omega _X=\Omega ^n_X$
is the canonical sheaf of differential forms of the highest degree.
For general $X$, $\omega ^\bcdot _X$ is known to be in $\dD (X)$.

For a closed embedding $i:X\to Y$ of $X$ into a smooth
$Y$ we have:
\begin{equation}\label{omxy}
\omega_X^{\bcdot}= i^!\omega_Y^{\bcdot}=i^!\omega_Y[{\rm dim}Y].
\end{equation}

The category $\dD (X)$ possesses a contravariant involutive exact functor
$D=D_X$.
%, the derived functor of local homomorphisms into $\omega^{\bcdot }_X$.
For an ${\cal F}\in \dD (X)$, it is defined by:
\begin{equation}
D{\cal F}=D_X{\cal F}={\cal H}om_X({\cal F}, \omega^{\bcdot}_X).
\end{equation}

$D$ is compatible with triangulated structures, i.e.
%it has the properties that
%\begin{itemize}
%\item[i)]
$DT$ is naturally isomorphic to $T^{-1}D$,
%\item[ii)]
and $D$ takes an exact triangle $A\to B\to C$
to the exact triangle
$DC\to DB\to DA$.
%\end{itemize}

If $j:U\to X$ is an open embedding, then
the dualizing complex over $U$ is obtained by restricting of the one on $X$:
\begin{equation}\label{jom}
\omega^{\bcdot}_U=j^*\omega^{\bcdot}_X.
\end{equation}
It follows that $D$ is local, meaning that, given an open embedding as above, we have:
\begin{equation}\label{jdu}
j^*D_X=D_Uj^*.
\end{equation}

We use notation $D^{(i)}{\cal F} := {\cal H}^i(D{\cal F})$
for the cohomology sheaves of $D{\cal F}$. $D$ preserves $\dD (X)$
and $D\cdot D$ is naturally isomorphic to the identity functor (\cite {RD}).

Using $D$, we can define another \tr e
${\tilde {\cal T}}=({\tilde \Dlz }, {\tilde \Dgo })$ on $\dD (X)$ by requiring:
\begin{equation}
{\cal F}\in {\tilde \Dlz }\ \Longleftrightarrow \ D{\cal F}\in {\cal D}^{\ge 0},
\end{equation}
\begin{equation}
{\cal F}\in {\tilde {\cal D}}^{\ge 0}\ \Longleftrightarrow \ D{\cal F}\in \Dlz.
\end{equation}
As $D$ is an anti-equivalence, ${\tilde {\cal T}}$ is a \tr e. We call it the
{\em dual \tr e}. Note that the dualizing complex is a pure object for the dual \tr e.

\subsection{A consistent pair of chains and perverse sheaves}

We denote by $\tau _{\le k}$, $\tau _{\ge k}$ (resp. ${\tilde \tau _{\le k}}$, ${\tilde \tau _{\ge k}}$)
the truncation functors for ${\cal T}$ (resp. ${\tilde {\cal T}}$).

The definition of  ${\tilde {\cal T}}$ implies:
\begin{equation}\label{dtld}
{\tilde \tau _{\le r}}=D\tau_{\ge -r}D,\ \ {\tilde \tau _{\ge r}}=D\tau_{\le -r}D.
\end{equation}

\begin{PROP}\label{ltc}
Let $X$ be a scheme of finite type over $k$. Suppose that $G\in {\cal D}^{\ge 0}(X)$,
then, for any $r\in\zz $, $D(\tau _{\le r}DG)\in {\cal D}^{\ge 0}(X)$.
\end{PROP}
\begin{PROOF}{.} In view of (\ref{jdu}), the statement of the proposition is local. Therefore, we may assume that
$X$ is affine and (being of finite type) embeddable into  a smooth variety over $k$.
Fix such a closed embedding $i:X\to Y$ of $X$ into a smooth variety $Y$ of dimension $l$. In view
of (\ref{adish}) and (\ref{omxy}), for any $G\in {\cal D}(X)$, we have:
$$
i_*DG=i_*{\cal H}om_X(G,\omega ^{\bcdot}_X)=i_*{\cal H}om_X(G,i^!\omega ^{\bcdot}_Y)={\cal H}om_Y(i_*G,\omega _Y[l]).
$$
Since $i_*$ is an exact functor with respect to the standard \tr e, then
\begin{equation}\label{idk}
i_*D^{(k)}G=\mathbb{R}^ki_*DG=\mathbb{R}^k{\cal H}om_Y(i_*G,\omega_Y[l]).
\end{equation}

Let $G\in {\cal D}^{\ge 0}(X)$. Then $\mathbb{R}^ki_*G$ are zero for $k<0$.

Consider the spectral sequence with  $E_2^{mj}={\cal E}xt^{m+l}_Y(\mathbb{R}^{-j}i_*G,\omega_Y)$
that converges to $\mathbb{R}^{m+j}{\cal H}om_Y(i_*G,\omega _Y[l])$.

Since $\omega _Y$ is locally free, then, for any coherent sheaf
${\cal F}$ on $Y$, ${\cal E}xt_Y^k({\cal F},\omega _Y)$ has
support of codimension
 $\ge k$ in $Y$ (cf. \cite{OSS}, p.142). By
convention, the support
of the zero sheaf, which is the empty set, is of infinite codimension.

As $E_2^{mj}=0$, for $j>0$, it follows from the spectral sequence that the
support of $\mathbb{R}^s{\cal H}om_Y(i_*G, \omega _Y[l])$
is of codimension at least $s+l$. Now, formula (\ref{idk}) implies that
the same restriction on codimension of support is verified for $i_*D^{(s)}G$.

As the cohomology sheaves of the truncation $\tau _{\le r}i_*DG$ either
coincide
with those of $i_*DG$ or equal zero, the codimension of the support
of ${\cal H}^s(\tau _{\le r}i_*DG)$
in $Y$ is at least $s+l$ for any $r\in\mathbb{Z}$.

Now consider the spectral sequence with $E_2^{js}=D^j{\cal H}^{-s}(\tau _{\le r}DG)$, which converges to the cohomology
sheaves ${\cal H}^{j+s}(D(\tau _{\le r}DG))$.

By (\ref{idk}) and in view of exactness of $i_*$ we obtain:
$$
i_*D^j{\cal H}^{-s}(\tau_{\le r}DG)=\mathbb{R}^j{\cal H}om_Y(i_*{\cal H}^{-s}(\tau _{\le r}DG), \omega_Y[l])={\cal E}xt_Y^{j+l}({\cal H}^{-s}(\tau_{\le r}i_*DG),\omega_Y).
$$
Using the fact
that ${\cal E}xt_Y^k({\cal F},{\cal E})=0$ for any locally free sheaf ${\cal E}$
and for any coherent sheaf ${\cal F}$ with codimension of the support
greater than $k$, we find that the right hand side is trivial for $j<-s$.

Then the spectral sequence implies
$$
i_*{\cal H}^{j}(D(\tau _{\le r}DG))=0.
$$
The statement of the proposition follows.
\end{PROOF}

Now consider the set of \tr es in $\dD (X)$ as a set with consistencies where consistency
binary relations are defined as in section \ref{pairs}.
The sequences $\tT [r]$ and ${\tilde {\tT }}[r]$ can be viewed as chains in this partially ordered set.
An ordered pair of \tr es is said to be {\em $T$-consistent} if the pair of chains obtained from these two \tr es by applying translation functor is consistent.

\begin{THM}\label{tccoh}
The pair $({\cal T}, {\tilde {\cal T}})$ is a $T$-consistent pair of \tr es.
\end{THM}
\begin{PROOF}{.}
We have to show that, for any $r\in\mathbb{Z}$, the pair $({\cal T}[r],{\tilde {\cal T}})$ is lower and upper consistent.

If, for $F\in {\tilde {\cal D}}^{\le 0}$, we denote $G=DF$, then $G\in {\cal D}^{\ge 0}$. By the proposition $D\tau_{\le r}DG=D\tau_{\le r}F$ belongs to ${\cal D}^{\ge 0}$.
Hence $\tau_{\le r}{\tilde \Dlz }\subset {\tilde \Dlz}$, which is the condition of lower consistency.

Further, in view of (\ref{dtld}), we have:
$$
{\tilde \tau}_{\ge r}{\cal D}^{\ge 0} =D\tau_{\le -r}D{\cal D}^{\ge 0}\subset {\cal D}^{\ge 0},
$$
which yields the upper consistency.
\end{PROOF}

Due to theorem \ref{tccoh}, pair of chains $\tT [r]$ and ${\tilde {\tT }}[r]$ is consistent.
Hence, by theorem \ref{latti} we obtain a distributive lattice of \tr es.
The \tr es that occur in this lattice are in fact coherent versions of
perverse \tr es (cf. \cite{Bez}).

Now recall that for \tr es there is the naive intersection (resp. union) defined by
intersecting the subcategories $\Dlz$ (resp. $\Dgo$).

\begin{PROP}
Any (multiple) intersection and union of \tr es in the obtained lattice
 is naive.
\end{PROP}
\begin{PROOF}{.}
Intersection or union in consistent pairs is manifestly naive. As follows from
 the remark after lemma \ref{strass}, the situation described in that lemma is the only place where we should check naivety of the union and intersection.
%in the course of the proof of theorem \ref{} was lemma \ref{}
So the fact comes with the following reformulation of lemma \ref{strass}
for \tr es:
if $a+b$ exists and is naive, then the existence and naivety for one side of the equation (\ref{abcabc}) implies existence and naivety for the other side and the equation itself.
This statement is again obvious.
\end{PROOF}

%Let us describe perverse \tr es in more familiar terms.

%%%%%%%%%%%%%%%%%%%%%%%%%%%%%%%%%%%%%%%%%

\end{document}